\documentclass[a4paper,10pt]{article}

\title{\LARGE  \bf Homoclinic Invariants of Ergodic Actions }
\author{Valery V. Ryzhikov}
\date{}

\begin{document}
\Large
\maketitle
\begin{abstract}{\large We consider a family of homoclinic groups and Gordin's 
type invariants of measure-preserving actions, state 
their connections with factors, full groups, ranks, rigidity, 
multiple mixing and   realize such invariants in the class of
  Gaussian and Poisson suspensions.}
\end{abstract}
\section{Introduction}
In a topological group $ G $ there is the correspondence  between an element 
 $ T $ and its homoclinic group:
$$ H (T) = \{S: T^{- n} ST^n \to I, \ n \to \infty \}, $$
 where $ I $ is the neutral element. We consider the case where $ G $ is
the  group of all automorphisms of a Lebesgue space. 
 It is not hard to see that all Bernoulli actions 
have ergodic homoclinic groups 
(the same is probably true for all K-automorphisms).
King, answering  Gordin's question, built in [2] a zero entropy transformation   $ T $ of 
a probability space
 with   ergodic $ H (T)$, we say  \it a transformation $ T $  with
 Gordin's G-property. \rm 
This invariant implies the  mixing  [1], furthermore, 
it implies the  mixing of all orders [8].

All mixing Gaussian and Poisson suspensions  (see [3] for definitions),
 have G-property (the proof  for Poisson suspensions see in  [8], 
as for Gaussian actions, it is an exercise). 
 
We know that all group actions without  multiple mixing property
and the horocycle flows 
do not possess  Gordin's property 
(the latter follows  from Ratner's results [6]). 

Let us define a family  of homoclinic nature groups.
 The weak homoclinic group is defined as 
$$ WH (T) = \{S\in G: \frac {1} {N} \sum_{n = 1}^{N}T^{- n} ST^n \to I, 
\ N \to \infty \}. $$

 For $ P $, an infinite subset of integers, we define the group
$$ H_P (T) = \{S\in G \ : \ T^{- n} ST^n \to I, \ n \in P, \ n \to \infty \}. $$

One presumes here the strong operator convergence, associating the   transformations  with corresponding 
 operators in $L_2(X,\mu)$.

\section{Results}
A transformation $ T $ is prime by defenition as  
it possesses only trivial invariant $\sigma$-algebras.
Recall that the full group $[S]$  of the automorphism  $S$ is defined as 
the collection of all  automorphisms $R$ of the Lebesgue space $(X,\mu)$
such that for all $x\in X$  one has $R(x)=S^{n(x)}$. 
 The full group $[\{S_g\}]$ of
a family $\{S_g\}$  is the group generated by all groups $[S_g]$.

\vspace {5mm}
\bf Theorem 1. \it If a transformation $ T $ is prime,
 then either the group  $H_P (T) $ is ergodic or $ H_P (T) = \{I \} $. 

  For all $P$ the group $ H_P (T) $ is full : $H_P (T) = 
[H_P (T)]$.

\bf Corollary. \it If the group $ WH (T) $ ($H (T) $,$ H_P (T) $) is ergodic,
then it includes representatives of all conjugacy classes of ergodic transformations.
King's homoclinic group from [2]  contains isomorphic copies of all ergodic transformations. \rm

\vspace {5mm}
\bf Theorem 2. \it Weakly mixing  Gaussian and Poisson suspensions
 have  ergodic weak homoclinic group.\rm
\vspace {5mm}

 Let $T$ be an automorphism of  a standard Lebesgue  space $(X,\mu)$, $\mu(X)=1$.
Let for a number  $\beta>0$
there exist a  sequence $\xi_j$  of  partitions of $X$  in the form
$$ \xi_j=\{ B_j,  TB_j,    T^2B_j,    \dots,  T^{h_j-1}B_j,C_j^1,\dots, C_j^{m_j}\dots\},$$
and  
any measurable set can be approximated by 
$\xi_j$-measurable ones as $j\to\infty$,  and 
 $\mu(U_j)\to \beta$, where 
$U_j=\bigsqcup_{0\leq k<h_j}T^kB_j$. 
\it
The local rank  $\beta(T)$  is defined as   maximal  number  $\beta$ 
for which the automorphism $T$ possesses a corresponding sequence of approximating partitions. 
 An automorphism $T$ is said to be of   rank $ 1$, if  $\beta(T)=1$.\rm
\vspace {5mm}

\bf Theorem 3. \it If for an automorphism $ T $ of rank 1 the  group $WH (T) $ is infinite,
then $T$ is rigid: for some sequence $ n_i \to \infty $ it
holds $ T^{n_i} \to I $. 

 If a mixing  transformation  $ R $ is of rank 1, then
$H (R) = WH (R) = \{I \} $. 

If an automorphism $ T $ is ergodic, $\beta(T)=\beta >0$  and the  group $WH (T) $ is infinite,
then there is a sequence $ n_i \to \infty $ such that   $ T^{n_i} \to \beta I + (1-\beta)M $
for some Markov operator $M$ (we say:   $T$ is partially rigid). 

\rm

\vspace {5mm}

Theorems 2,3 imply  a generalization of Parreau-Roy's result: rank-one Poisson suspension 
must be rigid [5] (for now there are no examples of such kind).
 \vspace {5mm}

\bf Theorem 4. \it For any two infinite families of integers
there  are some   subsets $Q$, $P$ of these families, respectively, such that for some weakly mixing 
Poisson 
(Gaussian) 
suspensions $T$,$T'$   the groups
$ H_P (T)$, $H_Q (T') $ are ergodic and the groups $ H_P (T')$, $H_Q (T)$ are trivial. \rm
\vspace {5mm} 

The homoclinic approach gives  new proofs of 
the multiple mixing for Poisson suspensions (Roy [7]) and Gaussian actions (Leonov [4]) as well 
as the proof of  
 weak multiple mixing for weakly  mixing  Gaussian and Poisson suspensions.

\vspace {5mm}
\bf Theorem 5. \it Suppose that an automorphism $ T$
satisfies the properties
$ H_P (T) $ is ergodic, and 
for  sequences $ m_i^1, \dots, m_i^k \in P,  \  m_i^1, m^k_i \to\infty $ 
the  convergence
$$ \mu (T^{m_i^1} B_1 \cap \dots \cap T^{m_i^k} B_k) \to \mu (B_1) \dots \mu (B_k) $$
holds for any measurable sets $ B_1, \dots, B_k $.

Then for any measurable sets $ B, B_1, \dots, B_k $ we have 
$$ \mu (B \cap T^{m_i^1} B_1 \cap \dots \cap T^{m_i^k} B_k) 
\to \mu (B) \mu (B_1) \dots \mu (B_k). $$
\rm
\section{Proofs}

{\bf Proof of Theorem 1.} It is not hard to see that  the algebra of the fixed sets with respect to a homoclinic group 
$H_P (T) $
is a factor of $T$. Indeed, if  
  $T^{-n_i}ST^{n_i}\to I$, then
$$T^{-n_i-1}ST^{n_i+1}\to I,\ \ T^{-n_i}T^{-1}STT^{n_i}\to I.$$
Thus, from $A=T^{-1}STA$  we obtain  $TA=STA$, hence, we have the first part of Theorem 1.

The second part of Theorem 1  is an exercise as well.  

\vspace {5mm}
{\bf   Proof of Theorem 2.}  Let a  weakly mixing   Poisson suspension $T_*$ be induced by infinite transformation $T$.
 The latter has  the  ergodic weak homoclinic group. It follows from
 the fact that all finite measure supports infinite transformations $F$  are weakly homoclinic:
$$  \frac{1} {N} \sum_{n = 1}^{N}T^{- n} FT^n \to I, 
\ N \to \infty. $$
 Analogous arguments for   mixing transformations see in [8].

For Gaussian actions our method is  similar:  we use the fact that  the group $FO$  of all 
"finite dimension" orthogonal operators are  dense in the group of 
all orthogonal operators on ${\it l}_2$. The term "finite dimension" orthogonal operator means that
this operator is in the form $U\oplus I$, where $U$ acts on a finite dimension space. 
For any weakly mixing Gaussian transformation $T$
the Gaussian image of the group $FO$ will be an ergodic subgroup  of the group  $WH(T)$.

\vspace {5mm}
The proof of Theorem 3  requires  the methods of [9].

\vspace {5mm}
{\bf Proof of Theorem 4.} It is an exercise to construct infinite transformations $T$ 
 mixing along some  subset $P$
and  rigid along some  subset $Q$. Then we apply the corresponding Poisson and Gaussian constructions.
\vspace {5mm}

{\bf Proof of Theorem 5.}
Let  for a measure $\nu$ we have 
$$ \mu (B \cap T^{m_i^1} B_1 \cap \dots \cap T^{m_i^k} B_k) 
\to \nu (B \times B_1\times\dots \times B_k) $$
for any measurable sets $B,B_1, \dots B_k. $ 

The projection of the measure  $\nu$ into the cube $X_1\times\dots \times X_k $,
where $X$ denotes the phase space of our transformation, is $\mu^k$.

The measure $\nu$ is invariant with respect to $S\times Id\times Id\dots \times Id $
for all $S\in  H_P (T)$. 
Indeed,  
$$ \mu (B \cap T^{m_i^1} B_1 \cap \dots \cap T^{m_i^k} B_k) $$
$$=\mu (SB \cap ST^{m_i^1} B_1 \cap \dots \cap ST^{m_i^k} B_k)$$
$$=\mu (SB \cap T^{m_i^1}T^{-m_i^1}ST^{m_i^1} B_1 
\cap \dots \cap T^{m_i^k}T^{-m_i^k}ST^{m_i^k} B_k). $$
But
$$\mu (SB \cap T^{m_i^1}T^{-m_i^1}ST^{m_i^1} B_1 
\cap \dots \cap T^{m_i^k}T^{-m_i^k}ST^{m_i^k} B_k)- \mu (SB \cap T^{m_i^1} B_1 
\cap \dots \cap T^{m_i^k}B_k)\to 0. $$
From this we get
$$ \mu (B \cap T^{m_i^1} B_1 \cap \dots \cap T^{m_i^k} B_k)\to \nu (SB \times B_1\times\dots \times B_k),$$
hence, for all  $S\in  H_P (T)$
$$\nu (B \times B_1\times\dots \times B_k)=\nu (SB \times B_1\times\dots \times B_k).$$
 The ergodicity of the group $H_P (T)$ 
implies now $$\nu=\mu\times \mu^{k}=\mu^{k+1}.$$ 
Thus,
$$ \mu (B \cap T^{m_i^1} B_1 \cap \dots \cap T^{m_i^k} B_k) 
\to \mu (B) \mu (B_1) \dots \mu (B_k). $$

\section{Questions}

 \it 

Let $T$ be weakly mixing rank-one transformation.
What  can we  say about homoclinic groups $WH(T)$, $H_P(T)$? \rm

We  note that $$\bigcup_{P, \ d(P)=1} H_P(T)\subset WH(T),$$ where $d(P)$
denotes the density of the set  $P$.
\it

 Could  $H_P(T)$ be ergodic for mixing rank-one transformation $T$,  or, more generally,
for a transformation with  minimal self-joinings? \rm

From Ageev's results we know that for  generic actions of the group
$(S,T: T^{-1}S^2 T=S)$ the generator $S$ is  rank-one.    
 In connection with the above question it is naturally
to ask: could such $T$ be a rank-one automorphism?  For all  known models the corresponding 
element $T$ is out of rank one.  In  generic case    for some infinite set $P\subset \{2,4,8,\dots\}$
 the group $H_P(T)$ is ergodic.

Only for transformations  $R$  with  discrete spectrum we know that 
  $H_P(R)=\{I\}$ for all infinite sets $P$. \it

Are there  weakly mixing  transformations  with the same property? \rm


E-mail: vryzh@mail.ru

\end{document}